\documentclass[a4paper,10pt]{article}

\usepackage[T1]{fontenc} 
\usepackage{amsmath,amsthm}
\usepackage[english]{babel}
\usepackage[dvips]{graphicx}
\usepackage{amsfonts,amssymb}
\usepackage{geometry}
\usepackage{hyperref}
\usepackage{stmaryrd}
\usepackage{fancyhdr}
\usepackage{url}

\newtheorem{prop}{Proposition}[section]
\newtheorem{LM}{Lemma}[section]
\newtheorem{thm}{Theorem}[section]
\newtheorem{df}{Definition}[section]
\newtheorem{cor}{Corollary}[section]
\newtheorem{conj}{Conjecture}[section]
\newtheoremstyle{pourlesremarques}{\topsep}{\topsep}{\normalfont}{}{\bfseries}{.}{ }{}
\theoremstyle{pourlesremarques}
\newtheorem{rem}{Remark}[section]
\newtheorem{concl}{Conclusion}[section]

\title {\textbf{Distinguished representations and exceptional poles of the Asai-L-function}}
\author{Nadir MATRINGE\footnote{Nadir Matringe, IMJ, 2 place Jussieu, F-75251, Paris Cedex 05. Email: matringe@math.jussieu.fr.}}

\begin{document}

\maketitle

\begin{abstract}
 Let $K/F$ be a quadratic extension of p-adic fields. We show that a generic irreducible representation of $GL(n,K)$ is distinguished if and ony if its Rankin-Selberg Asai L-function has an exceptional pole at zero. We use this result to compute Asai L-functions of ordinary irreducible representations of $GL(2,K)$. In the appendix, we describe supercuspidal dihedral representations of $GL(2,K)$ in terms of Langlands parameter.

\end{abstract}

\section*{Introduction}

For $K/F$ a quadratic extension of local fields, let $\sigma$ be the
 conjugation relative to this extension, and $\eta _{K/F}$ be the
character of $F^*$ whose kernel is the set of norms from $K^*$. The conjugation $\sigma$ extends naturally to an automorphism of $GL(n,K)$, which we also denote by $\sigma$. If $\pi$ is a representation of $GL(n,K)$, we denote by $\pi^{\sigma}$ the representation $g\mapsto \pi(\sigma(g))$.\\
If $\pi$ is a smooth irreducible representation of $GL(n,K)$, and $\chi$ a character of $F^*$, the
 dimension of the space of linear forms on its space, which transform by $\chi$ under $GL(n,F)$ (with respect to the action $[(L,g)\mapsto L\circ\pi(g)]$), is known to be at
 most one (Proposition 11, \cite{F1}).
One says that $\pi$ is $\chi$-distinguished if this dimension is one, and says that $\pi$ is distinguished if it is $1$-distinguished.\\
Jacquet conjectured two results about distinguished representations of $GL(n,K)$.
Let $\pi$ be a smooth irreducible representation of $GL(n,K)$ and
 $\pi^{\vee}$ its contragredient.
The first conjecture states that it is equivalent for $\pi$ with central character trivial on $F^*$ to
be isomorphic to ${\pi^{\vee}}^{\sigma}$ and for $\pi$ to be distinguished or $\eta _{K/F}$-distinguished.
 In \cite{K}, Kable proved it for discrete series representations, using Asai $L$-functions.\\ 
The second conjecture, which is proved in \cite{K}, states that if  $\pi$ is a
 discrete series representation, then it cannot be distinguished and $\eta _{K/F}$-distinguished at the same time.\\
One of the key points in Kable's proof is that if a discrete
 series representation of $GL(n,K)$ is such that its Asai $L$-function has a pole at zero, then it is distinguished, Theorem 1.4 of \cite{AKT} shows that it is actually an equivalence.
This theorem actually shows that Asai $L$-functions of tempered distinguished representations
 admit a pole at zero.\\ 
In this article, using a result of Youngbin Ok which states that for a
 distinguished representation, linear forms invariant under the
 affine subgroup of $GL(n,F)$ are actually $GL(n,F)$-invariant (which
 generalises Corollary 1.2 of \cite{AKT}), we prove in Theorem
 \ref{pole2} that a generic representation is distinguished if and only if its Asai $L$-function admits an exceptional pole at zero. A pole at zero is always exceptional for Asai $L$-functions of discrete
 series representations (see explanation before Proposition \ref{K}). As a first application, we give in Proposition \ref{Lcusp} a formula for Asai $L$-functions of supercuspidal representations of $GL(n,K)$.\\
There are actually three different ways to define Asai $L$-functions: one via the local Langlands correspondence and in terms of Langlands parameters
denoted by $L_{W}(\pi,s)$, the one we use via the theory of Rankin-Selberg integrals denoted
by $L_{As}(\pi,s)$, and the Langlands-Shahidi method applied to a suitable unitary group, denoted by $L_{As,2}(\pi,s)$ (see \cite{AR}). It is expected that the above three $L$-functions
are equal.\\
For a discrete series representation $\pi$, it is shown in \cite{He} that $L_{W}(\pi,s)=L_{As,2}(\pi,s)$, and in \cite{AR} that $L_{As}(\pi,s)=L_{As,2}(\pi,s)$, both proofs using global methods.\\  
As a second application of our principal result, we show (by local methods) in Theorem \ref{egal} of Section \ref{gl(2)} that for an ordinary representation (i.e. corresponding through Langlands correspondence to an imprimitive 2 dimensional representation of the Weil-Deligne group) $\pi$ of $GL(2,K)$, we have $L_{W}(\pi,s)=L_{As}(\pi,s)$ . We recall that for odd residual characteristic, every smooth irreducible infinite-dimensional representation of $GL(2,K)$ is ordinary.\\
In the appendix (Section \ref{appendix}), we describe in Theorem \ref{distcusp} distinguished dihedral supercuspidal representations, this description is used in Section \ref{gl(2)} for the computation of $L_{As}$ for such representations.

\section{Preliminaries}

Let $E_1$ be a field, and $E_2$ a finite galois extension of $E_1$, we denote by $Gal(E_2/E_1)$ the Galois group of $E_2$ over $E_1$, and we denote by $Tr_{E_2/E_1}$ (respectively $N_{E_2/E_1}$) the trace (respectively the norm) function from $E_2$ to $E_1$. If $E_2$ is quadratic over $E_1$, we denote by $\sigma_{E_2/E_1}$ the non trivial element of $Gal(E_2/E_1)$.\\
In the rest of this paper, the letter $F$ will always designate a non archimedean local field of characteristic zero in a fixed algebraic closure $\bar{F}$, and the letter $K$ a quadratic extension of $F$ in $\bar{F}$. We denote by $q_F$ and $q_K$ the cardinality of their residual fields, $R_K$ and $R_F$ their integer rings, $P_K$ and $P_F$ the maximal ideals of $R_K$ and $R_F$, and $U_K$ and $U_F$ their unit groups. We also denote by $v_K$ and $v_F$ the respective normalized valuations, and $| \ |_K$ and $|\ |_F$ the respective absolute values. We fix an element $\delta$ of $K-F$ such that $\delta ^2 \in F$, hence $K=F(\delta)$.\\
Let $\psi$ be a non trivial character of $K$ trivial on $F$, it is of the form $x \mapsto \psi'\circ Tr_{K/F}(\delta x)$ for some non trivial character $\psi'$ of $F$.\\
Whenever $G$ is an algebraic group defined over $F$, we denote by $G(K)$ its
 $K$-points and $G(F)$ its $F$-points.
The group $GL(n)$ is denoted by $G_n$, its standard maximal unipotent
 subgroup is denoted by $N_n$.\\
If $\pi$ is a representation of a group, we also denote by $\pi$ its isomorphism class.
Let $\mu$ be a character of $F^*$, we say that a representation $\pi$ of $G_n(K)$ is $\mu$-distinguished if it admits on its space $V_{\pi}$ a linear form $L$, which verifies the following: for $v$ in $V$ and $h$ in $G_n(K)$, then $L(\pi(h)v)=\mu(det(h))L(v)$. If $\mu=1$, we say that $\pi$ is distinguished.\\
We denote by $K_n(F)$ the maximal compact subgroup $G_n(R_F)$ of $G_n(F)$, and for $r\geq 1$, we denote by $K_{n,r}(F)$, the congruence subgroup $I_n+ M_n(P_F^r)$.\\
The character $\psi$ defines a character of $N_n(K)$ that we still denote by
 $\psi$, given by $\psi(n)=\psi(\sum_{i=1}^{n-1} n_{i,i+1})$.\\

We now recall standard results from \cite{F2}.\\
Let $\pi$ be a generic smooth irreducible representation of $G_n(K)$,
 we denote by $\pi^{\vee}$ its smooth contragredient, and $c_{\pi}$ its
 central
 character.\\
We denote by $D(F^n)$ the space of smooth functions with compact support on
 $F^n$, and $D_0(F^n)$ the subspace of $D(F^n)$ of functions vanishing
 at zero.
We denote by $\rho$ the natural action of $G_n(F)$ on $D(F^n)$, given by
 $\rho(g)\phi (x_1,\dots,x_n)= \phi ((x_1,\dots,x_n)g)$, and we denote by $\eta$ the row vector $(0,\dots,0,1)$ of length $n$.\\
If $W$ belongs to the Whittaker model $W(\pi,\psi)$ of $\pi$, and
 $\phi$ belongs to $D(F^n)$, the following integral converges for $s$
 of
 real part large enough:

$$\int_{N_n(F)\backslash G_n(F)} W(g)\phi(\eta g){|det(g)|_F}^s dg.$$

This integral as a function of $s$ has a meromorphic extension to
 $\mathbb{C}$ which we denote by $\Psi(W,\phi,s)$.\\
For $s$ of real part large enough, the function $\Psi(W,\phi,s)$ is a
 rational function in $q_F^{-s}$, which actually has a Laurent series
 development.\\
The $\mathbb{C}$-vector space generated by these functions is in fact
 a fractional ideal $I(\pi)$ of $\mathbb{C}[q_F^{-s},q_F^s]$.
This ideal $I(\pi)$ is principal, and has a unique generator of the form
 $1/P(q_F^{-s})$, where $P$ is a polynomial with $P(0)=1$.

\begin{df}
We denote by $L_{As}(\pi,s)$ the generator of $I(\pi)$ defined just above,
 and call it the Asai $L$-function of $\pi$.
\end{df}

\begin{rem}\label{Eulerpole} If $P$ belongs to $\mathbb{C}[X]$ and has constant term equal to one, then the function of the complex variable $L_P:s\mapsto 1/P(q_F^{-s})$ is called an Euler factor. It is a meromorphic function on $\mathbb{C}$ and admits $(2i\pi/ln(q_F))\mathbb{Z}$ as a period subgroup. Hence if $s_0$ is a pole of $L_{P}$, the elements $s_0 + (2i\pi/ln(q_F))\mathbb{Z}$ are also poles of $L_P$, with same multiplicities, we identify $s_0$ and $s_0 + (2i\pi/ln(q_F))\mathbb{Z}$ when we talk about poles. A pole $s_0$ then corresponds to a root $\alpha_0$ of $P$ by the formula $q^{-s_0}=\alpha_0$, its multiplicity in $L_P$ equal to the multiplicity of $\alpha_0$ in $P$.\end{rem}

Let $w_n$ be the matrix of $G_n(\mathbb{Z})$ with ones on the antidiagonal, and zeroes elsewhere.
For $W$ in $W(\pi,\psi)$, we denote by $\tilde{W}$ the function $g\mapsto
 W(w_n {^t{g}^{-1}})$ which belongs to $W(\pi^{\vee},{\psi}^{-1})$,
 and
 we
 denote by $\widehat{\phi}$ the Fourier transform (with respect to $\psi'$ and its associate autodual Haar measure)
 of $\phi$ in $D(F^n)$.\\

\begin{thm}{(\textbf{Functional equation})}{(Th. of \cite{F2})}

There exists an epsilon factor $\epsilon _{As}(\pi,s,\psi)$ which
 is, up to scalar, a (maybe negative) power of $q^s$, such that the
 following functional equation is satisfied for any $W$ in
 $W(\pi,\psi)$ and
 any $\phi$ in $D(F^n)$:

$$\Psi(\tilde{W},\widehat{\phi},1-s)/L_{As}(\pi^{\vee},1-s) =
 c_{\pi}(-1)^{n-1} \epsilon
 _{As}(\pi,s,\psi)\Psi(W,\phi,s)/L_{As}(\pi,s).$$

\end{thm}

We finally recall the following, which will be crucial in the demonstration of Theorem \ref{pole2}.

\begin{prop}{(\cite{Ok}, Theorem $3.1.2$)} \label{bernok}
Let $\pi$ be an irreducible distinguished representation of $G_n(K)$, if $L$ is a $P_n(F)$-invariant linear form on the space of $\pi$, then it is actually $G_n(F)$-invariant.
\end{prop}

\begin{proof}[Sketch of the proof] We note $V$ the space of $\pi$, and $\tilde{V}$ that of $\pi^{\vee}$. As the representation $\pi^{\vee}$ is isomorphic to $g\mapsto \pi((g^t)^{-1})$, it is also distinguished. Let $L$ be a $P_n(F)$-invariant linar form on the space $V$ and $\tilde{L}$ a $G_n(F)$-invariant linar form on the space $\tilde{V}$, the linear form $L\otimes \tilde{L}$ on $V \otimes \tilde{V}$ is $P_n(F)\times G_n(F)$-invariant. It is thus enough to prove that a linear form $B$ on $V \otimes \tilde{V}$ which is $P_n(F)\times G_n(F)$-invariant is $G_n(F)\times G_n(F)$-invariant.\\
Call $\lambda$ the (right) action by left translation and $\rho$ that by right translation of $G_n(K)$ on the space $C_c^{\infty}(G_n(K))$, it follows from Lemma p.73 of \cite{B} that there exists an injective morphism $I$ of $G_n(K)\times G_n(K)$-modules from $[(\pi \otimes \pi^{\vee})^*, (V \otimes \tilde{V})^*]$ to $[(\lambda \times \rho)^*, (C_c^{\infty}(G_n(K)))^*]$. The linear form $I(B)$ is an element of $(C_c^{\infty}(G_n(K)))^*$ which is $P_n(F)\times G_n(F)$-invariant. As $I$ is injective, the result will follow from the fact that an invariant distribution on $G_n(K)/G_n(F)$ which is invariant by left translation under $P_n(F)$ is actually $G_n(F)$-invariant. Identifying $G_n(K)/G_n(F)$ with the space $S$ of matrices $g$ of $G_n(K)$ verifying of $gg^{\sigma}=1$ (see \cite{S}, ch.10, prop.3), this statement is exactly the one of Lemma 5 of \cite{GJR}. \end{proof}

\section{Poles of the Asai $L$-function and distinguishedness}

Now suppose $L_{As}(\pi,s)$ has a pole at $s_0$, its order $d$ is the
 highest order pole of the family of functions of $I(\pi)$.\\
Then we have the following Laurent expansion at $s_0$:  

\begin{equation}\label{dvlaurent} \Psi(W,\phi,s)=
 B_{s_0}(W,\phi)/(q_F^s-q_F^{s_0})^d + \ smaller \ order \
 terms. \end{equation}

The residue $B_{s_0}(W,\phi)$ defines a non zero bilinear form on
 $W(\pi,\psi) \times D(F^n)$, satisfying the quasi-invariance:

$$ B_{s_0}(\pi(g)W,\rho(g)\phi)= |det(g)|_F^{- s_0} B_{s_0}(W,\phi).$$

Following \cite{CP} for the split case $K=F\times F$, we state the following definition:

\begin{df}
 
A pole of the Asai $L$-function $L_{As}(\pi,s)$ at $s_0$ is called exceptional if the
 associated bilinear form $B_{s_0}$ vanishes on $W(\pi,\psi) \times
 D_0(F^n)$.\\

\end{df}

As an immediate consequence, if $s_0$ is an exceptional pole of
 $L_{As}(\pi,s)$, then $B_{s_0}$ is of the form $B_{s_0}(W,\phi)=
 \lambda
 _{s_0}(W) \phi(0)$, where $\lambda _{s_0}$ is a non zero $|det(\ )|_F^{-
 s_0}$
 invariant linear form on $W(\pi,\psi)$.\\
Hence we have:

\begin{prop}
\label{pole1}

Let $\pi$ be a generic irreducible representation of $G_n(K)$, and
 suppose its Asai $L$-function has an exceptional pole at zero, then
 $\pi$
 is distinguished.

\end{prop}

We denote by $P_n(F)$ the affine subgroup of $G_n(F)$, given by matrices with last row equal to $\eta$.\\
For more convenience, we introduce a second $L$-function: for $W$ in $W(\pi,\psi)$, by standard arguments, the following integral
 is convergent for $Re(s)$ large, and defines a rational function in $q^{-s}$, which has a Laurent series development:

$$\int_{N_n(F)\backslash P_n(F)} W(p){|det(p)|_F}^{s} dp.$$

We denote by $\Psi_1(W,s)$ the corresponding Laurent series.
By standard arguments again, the vector space generated by the
functions $\Psi_1(W,s-1)$, for $W$ in $W(\pi,\psi)$, is a fractional ideal $I_1 (\pi)$ of
 $\mathbb{C}[q_F^{-s},q_F^s]$,
 which has a unique generator of the form $1/Q(q_F^{-s})$, where $Q$ is a
 polynomial with $Q(0)=1$. We denote by $L_1(\pi,s)$ this generator.

\begin{LM}\label{LM}(\cite{JPS} p. 393)

Let $W$ be in $W(\pi,\psi)$, one can choose $\phi$ with support small enough around $(0,\dots,0,1)$ such that $\Psi(W,\phi,s)=\Psi_1(W,s-1)$.
 
\end{LM}
\begin{proof} As we gave a reference, we only sketch the proof.
We first recall the following integration formula (cf. proof of the proposition in paragraph 4 of \cite{F}), for $Re(s)>>0$: 

\begin{equation} \label{sum1} \Psi(W,{\phi},s)=
 \int_{K_n(F)}
 \int_{N_n(F)\backslash P_n(F)} W(pk)|det(p)|_F^{s-1}dp \int_{F^*}
 \phi(\eta
 ak)c_{\pi}(a)|a|_F^{ns} d^*a dk.\end{equation}

Choosing $r$ large enough for $W$ to be right invariant under $K_{n,r}(F)$, we take $\phi$ a positive multiple of the characteristic function of $\eta K_{n,r}(F)$, and conclude from equation (\ref{sum1}).\end{proof} 

Hence we have the inclusion $I_1 (\pi) \subset I(\pi)$, which implies
 that $L_1(\pi,s)= L_{As}(\pi,s) R(q_F^s,q_F^{-s})$ for some $R$ in $\mathbb{C}[q_F^{-s},q_F^s]$.
But because $L_1$ and $L_{As}$ are both Euler factors, $R$ is actually
 just a polynomial in $q_F^{-s}$, with constant term equal to one.
Noting $L_{rad(ex)} (\pi,s)$ its inverse (which is an Euler factor), we have
 $L_{As}(\pi,s)=L_1(\pi,s)L_{rad(ex)} (\pi,s)$, we will say that $L_1$ divides
 $L_{As}$. The explanation for the notation $L_{rad(ex)}$ is given in Remark \ref{simplepoles}.\\

We now give a characterisation of exceptional poles:

\begin{prop}
 
A pole of $L_{As}(\pi,s)$ is exceptional if and only if it is a pole of the function $L_{rad(ex)} (\pi,s)$ defined just above.

\end{prop}

\begin{proof}
From equation (\ref{sum1}), it becomes clear that the vector space generated by the integrals $\Psi(W,\phi,s)$
 with $W$ in $W(\pi,\psi)$ and $\phi$ in $D_0(F^n)$, is contained in $I_1 (\pi)$, but because of Lemma \ref{LM}, those two vector spaces are equal.
Hence $L_1(\pi,s)$ is a generator of the
 ideal generated as a vector space by the functions $\Psi(W,\phi,s)$
 with
 $W$ in $W(\pi,\psi)$ and $\phi$ in $D_0(F^n)$.\\
From equation (\ref{dvlaurent}), if $s_0$ is an exceptional pole, a
 function $\Psi(W,\phi,s)$, with $\phi$ in $D_0(F^n)$, cannot have a pole of highest order at $s_0$, hence we have
 one
 implication.\\
Now if the order of the pole $s_0$ for $L_{As}(\pi,s)$ is stricly
 greater than the one of $L_1(\pi,s)$, then the first residual term corresponding to a pole of highest order of the Laurent development of any function $\Psi(W,\phi,s)$ with $\phi(0)=0$ must be zero, and zero is exceptional.\end{proof} 

Lemma \ref{LM} also implies:\\

\begin{prop}
 
The functional ${\Lambda}_{\pi,s}: W \mapsto
  \Psi_1(W,s-1)/L_{As}(\pi,s)$ defines a (maybe null) linear form on
 $W(\pi,\psi)$ which transforms
 by $|det(\ )|_F^{1-s}$ under the affine subgroup $P_n(F)$.\\
For fixed $W$ in $W(\pi,\psi)$, then $s \mapsto {\Lambda}_{\pi,s}(W)$ is a polynomial of
 $q_F^{-s}$.

\end{prop}

Now we are able to prove the converse of Proposition \ref{pole1}:

\begin{thm}\label{pole2} 
A generic irreducible representation $\pi$ of $G_n(K)$ is
 distinguished if and only if $L_{As}(s,\pi)$ admits an exceptional
 pole at zero.

\end{thm}

\begin{proof} We only need to prove that if $\pi$ is distinguished, then
 $L_{As}(s,\pi)$ admits an exceptional pole at zero, so we suppose
 $\pi$
 distinguished.\\
From equation (\ref{sum1}), for $Re(s)<<0$, and $\pi$ distinguished (so that $c_{\pi}$ has trivial restriction to $F^*$), one has:

\begin{equation} \label{sum} \Psi(\tilde{W},\widehat{\phi},1-s)=
 \int_{K_n(F)}
 \int_{N_n(F)\backslash P_n(F)} \tilde{W}(pk)|det(p)|_F^{-s}dp \int_{F^*}
 \widehat{\phi}(\eta
 ak)|a|_F^{n(1-s)} d^*a dk.\end{equation}

This implies that:

\begin{equation} \label{summ} \Psi(\tilde{W},\widehat{\phi},1-s)/L_{As}(\pi^{\vee},1-s)=\int_{K_n(F)}{\Lambda}_{\pi^{\vee},1-s}(\pi^{\vee}(k)\tilde{W})\int_{F^*}
 \widehat{\phi}(\eta ak)|a|_F^{n(1-s)} d^*a dk.\end{equation}

The second member of the equality is actually a finite sum: $\sum_i \lambda
 _i{\Lambda}_{\pi^{\vee},1-s}(\pi^{\vee}(k_i)\tilde{W})\int_{F^*}
 \widehat{\phi}(\eta ak_i)|a|_F^{n(1-s)} d^* a$, where the
 $\lambda_i$'s are positive constants and the $k_i$'s are elements of $K_n(F)$ independant of $s$.\\
Note that there exists a positive constant $\epsilon$, such that for
 $Re(s)<\epsilon$, the integral $\int_{F^*} \widehat{\phi}(\eta
 ak_i)|a|_F^{n(1-s)} d^* a$ is
 absolutely convergent, and defines a holomorphic function.
So we have an equality (equality \ref{summ}) of analytic functions (actually of polynomials in
 $q_F^{-s}$), hence it is true for all $s$ such that $Re(s)<\epsilon$.\\
For $s=0$, we get:

 $$\Psi(\tilde{W},\widehat{\phi},1)/L_{As}(\pi^{\vee},1)=\int_{K_n(F)}{\Lambda}_{\pi^{\vee},1}(\pi^{\vee}(k)\tilde{W})\int_{F^*} \widehat{\phi}(\eta
 ak)|a|_F^{n} d^*a dk.$$

But as $\pi$ is distinguished, so is $\pi^{\vee}$, and as
 ${\Lambda}_{\pi^{\vee},1}$ is a $P_n(F)$-invariant linear form on
 $W(\pi^{\vee},\psi^{-1})$, it follows from Propodsition
\ref{bernok} that it is
 actually $G_n(F)$-invariant.\\
Finally $$\Psi(\tilde{W},\widehat{\phi},1)/L_{As}(\pi^{\vee},1)=
 {\Lambda}_{\pi^{\vee},1}(\tilde{W})\int_{K_n(F)}\int_{F^*}
 \widehat{\phi}(\eta ak)|a|_F^{n} d^*a dk$$ which is equal to:\\
$${\Lambda}_{\pi^{\vee},1}(\tilde{W})\int_{P_n(F)\backslash
 G_n(F)}\widehat{\phi}(\eta g)d_{\mu}g$$ where $d_{\mu}$ is up to
 scalar the unique $|det(\ )|^{-1}$ invariant measure on
 $P_n(F)\backslash
 G_n(F)$.
 But as $$\int_{P_n(G)\backslash
 G_n(F)}\widehat{\phi}(\eta g)d_{\mu}g=
 \int_{F^n}\widehat{\phi}(x)dx=\phi(0),$$
we deduce from the functional equation that
 $\Psi(W,\phi,0)/L_{As}(\pi,0)=0$ whenever $\phi(0)=0$.\\
As one can choose $W$, and $\phi$ vanishing at zero, such that
  $\Psi(W,\phi,s)$ is the constant function equal to $1$ (see the proof of Theorem 1.4 in \cite{AKT}),
 hence $L_{As}(\pi,s)$ has a pole at zero, which must be
 exceptional. \end{proof} 

For a discrete series representation $\pi$, it follows from Lemma 2
 of \cite{K}, that the integrals of the form $$\int_{N_{n}(F)\backslash
 P_n(F)} W(p){|det(p)|_F}^{s-1} dp.$$ converge absolutely for
 $Re(s)>-\epsilon$ for some positive $\epsilon$, hence as functions of
 $s$, they
 cannot have a pole at zero.\\
This implies that $L_1(\pi,s)$ has no pole at zero, hence Theorem \ref{pole2} in this case gives:\\

\begin{prop}\label{K}(\cite{K}, Theorem 4)\\
 
A discrete series representation $\pi$ of $G_n(K)$ is
 distinguished if and only if $L_{As}(s,\pi)$ admits a pole at zero.

\end{prop}

Let $s_0$ be in $\mathbb{C}$.
We notice that if $\pi$ is a generic irreducible representation of $G_n(K)$, it is $|\ |_{F}^{-s_0}$-distinguished if and only if $\pi \otimes |\ |_{K}^{s_0 /2}$ is distinguished, but as $L_{As}(s,\pi \otimes |\ |_{K}^{s_0 /2})$ is equal to $L_{As}(s+s_0,\pi)$, Theorem \ref{pole2} becomes:

\begin{thm}\label{s-dist}

A generic irreducible representation $\pi$ of $G_n(K)$ is
 $|\ |_{F}^{-s_0}$-distinguished if and only if $L_{As}(s,\pi)$ admits an exceptional
 pole at $s_0$.
 
\end{thm}

\begin{rem}\label{simplepoles}

Let $P$ and $Q$ be two polynomials in $\mathbb{C}[X]$ with constant term 1, we say that the Euler factor $L_P (s)=1/P(q_F^{-s})$ divides $L_Q (s)=1/Q(q_F^{-s})$ if and only $P$ divides $Q$. We denote by $L_P \vee L_Q$ the Euler factor $1/(P\vee Q)(q_F^{-s})$, where the l.c.m $P\vee Q$ is chosen such that $(P\vee Q)(0)=1$. We define the g.c.d $L_P \wedge L_Q$ the same way.\\
It follows from equation (\ref{sum1}) that if ${c_{\pi}}_{|F^*}$ is ramified, then $L_{As}(\pi,s)=L_1(\pi,s)$.
It also follows from the same equation that if ${c_{\pi}}_{|F^*}=|\ |_F^{-s_1}$ for some $s_1$ in $\mathbb{C}$, then $L_{rad(ex)}(\pi,s)$ divides $1/(1-q_F^{s_1 -ns})$. Anyway, $L_{rad(ex)}(\pi,s)$ has simple poles.\\
Now we can explain the notation $L_{rad(ex)}$. We refer to \cite{CP} where the case $K=F\times F$ is treated. In fact, in the latter, $L_{ex}(\pi,s)$ is the function $1/P_{ex}(\pi,q_F^{-s})$, with $P_{ex}(\pi,q_F^{-s})=\prod_{s_i} (1-q_F^{s_i-s})^{d_i}$, where the $s_i$'s are the exceptional poles of $L_{As}(\pi,s)$ and the $d_i$'s their order in $L_{As}(\pi,s)$. Hence $L_{rad(ex)}(\pi,s)=1/P_{rad(ex)}(\pi,q_F^{-s})$, where $P_{rad(ex)}(\pi,X)$ is the unique generator with constant term equal to one, of the radical of the ideal generated by $P_{ex}(\pi,X)$ in $\mathbb{C}[X]$.\end{rem}

We proved:

\begin{prop}\label{Lrad(ex)}

Let $\pi$ be an irreducible generic representation of $G_n(K)$, the Euler factor $L_{rad(ex)}(\pi,s)$ has simple poles, it is therefore equal to $\prod 1/(1-q_F^{s_0-s})$ where the product is taken over the $q_F^{s_0}$'s such that $\pi$ is $|\ |_F^{-s_0}$-distinguished. 

\end{prop}

Suppose now that $\pi$ is supercuspidal, then the restriction to $P_n(K)$ of any $W$ in $W(\pi,\psi)$ has compact support modulo $N_n(K)$, hence $\Psi_1(W,s-1)$ is a polynomial in $q^{-s}$, and $L_1(\pi,s)$ is equal to $1$. Hence Proposition \ref{Lrad(ex)} becomes: 

\begin{prop}\label{Lcusp}
Let $\pi$ be an irreducible supercuspidal representation of $G_n(K)$, then $L_{As}(\pi,s)=\prod 1/(1-q^{s_0-s})$ where the product is taken over the $q^{s_0}$'s such that $\pi$ is $|\ |_F^{-s_0}$-distinguished.
\end{prop}

\section{Asai $L$-functions of $GL(2)$}\label{gl(2)}

\subsection{Asai $L$-functions for imprimitive Weil-Deligne
 representations of dimension 2}\label{Weil}

The aim of this paragraph is to compute $L_W(\rho,s)$ (see the introduction) when $\rho$ is an imprimitive two dimensional representation of the Weil-Deligne group of $K$.\\

We denote by $W_K$ (resp. $W_F$) the Weil group of $K$ (resp. $F$), $I_K$ (resp. $I_F$) the inertia subgroup of $W_K$ (resp. $W_F$), $W'_K$ (resp. $W'_F$) the group $W_K \times SL(2,\mathbb{C})$ (resp.
 $W_F
 \times SL(2,\mathbb{C})$) and $I'_K$ (resp. $I'_F$) the group $I_K \times SL(2,\mathbb{C})$ (resp. $I_F \times SL(2,\mathbb{C})$). We denote by $\phi_F$ a Frobenius element of $W_F$, and we also denote by $\phi_F'$ the element $(\phi_F,I_2)$ of $W'_F$.\\ 
We denote by $sp(n)$ the unique (up to isomorphism) complex irreducible representation of
 $SL(2,\mathbb{C})$ of dimension $n$.\\
If $\rho$ is a finite dimensional representation of $W'_K$, we denote by
 $M_{W'_K}^{W'_F}(\rho)$ the representation of $W'_F$ induced
 multiplicatively from $\rho$. We recall its definition:\\
If $V$ is the space of $\rho$, then the space of
 $M_{W'_K}^{W'_F}(\rho)$ is $V\otimes V$. Noting $\tau$ an element of $W_F -W_K$, and $\sigma$ the element
 $(\tau, I)$ of $W'_F$, we have:\\

$$ M_{W'_K}^{W'_F}(\rho)(h)(v_1 \otimes v_2) = \rho(h) v_1 \otimes
 \rho^{\sigma}(h) v_2$$ for $h$ in $W'_K$, $v_1$ and $v_2$ in $V$.

$$ M_{W'_K}^{W'_F}(\rho)(\sigma)(v_1 \otimes v_2) = \rho(\sigma ^2) v_2
 \otimes v_1$$ for $v_1$ and $v_2$ in $V$.\\

We refer to paragraph 7 of \cite{P} for definition and basic properties
 of multiplicative induction in the general case.

\begin{df}
The function $L_W(\rho,s)$ is by definition the usual $L$-function of the representation $M_{W'_K}^{W'_F}(\rho)$, i.e. $L_W(\rho,s)=L(M_{W'_K}^{W'_F}(\rho),s)$.
\end{df}

\begin{description}
 \item[i)] If $\rho$ is of the form $Ind_{W'_B}^{W'_K}(\omega)$ for
 some multiplicative character $\omega$ of a biquadratic extension $B$ of
 $F$, we denote by $K'$ and $K''$ the two other extensions between $F$ and
 $B$.
If we call $\sigma_1$ an element of $W'_K$ which is not in $W'_{K'}\cup W'_{K''}$ and $\sigma_3$ an element of $W'_{K''}$ which is not in $W'_K \cup W'_{K'}$, then $\sigma_2=\sigma_3\sigma_1$ is an element of $W'_{K'}$ which is not in $W'_K \cup W'_{K''}$.\\
The elements $(1,\sigma_1,\sigma_2,\sigma_3)$ are representatives of $W'_F/W'_B$, and $1$ and $\sigma_3$ are representatives of $W'_F/W'_K$.\\
 If one identifies $\omega$ with a character (still called $\omega$) of $B^*$, then $\omega^{\sigma_1}$ identifies with $\omega \circ \sigma_{B/K}$, $\omega^{\sigma_2}$ with $\omega \circ \sigma_{B/K'}$ and $\omega^{\sigma_3}$ with $\omega \circ \sigma_{B/K''}$. One then verifies that if $a$ belongs to $W_B$, one has:\\
$\bullet$ $Tr [ M_{W'_K}^{W'_F}(\rho)(a)] =Tr [Ind_{W'_{K'}}^{W'_F}(M_{W'_B}^{W'_{K'}}(\omega)) (a)]+Tr[Ind_{W'_{K''}}^{W'_F}(M_{W'_B}^{W'_{K''}}(\omega))(a)] = \omega\omega^{\sigma_2}+ \omega\omega^{\sigma_3}+ \omega^{\sigma_1}\omega^{\sigma_2}+ \omega^{\sigma_1}\omega^{\sigma_3}$.\\
$\bullet$ $Tr[M_{W'_K}^{W'_F}(\rho)(\sigma_1 a)]= Tr[Ind_{W'_{K'}}^{W'_F}(M_{W'_B}^{W'_{K'}}(\omega))(\sigma_1 a)]+Tr[Ind_{W'_{K''}}^{W'_F}(M_{W'_B}^{W'_{K''}}(\omega))(\sigma_1 a)]=0$.\\
$\bullet$ $Tr[M_{W'_K}^{W'_F}(\rho)(\sigma_2 a)]= Tr[Ind_{W'_{K'}}^{W'_F}(M_{W'_B}^{W'_{K'}}(\omega))(\sigma_2 a)]+Tr[ Ind_{W'_{K''}}^{W'_F}(M_{W'_B}^{W'_{K''}}(\omega))(\sigma_2 a)]= \omega(\sigma_2 a\sigma_2 a)+\omega^{\sigma_1}(\sigma_2 a\sigma_2 a)$.\\
$\bullet$ $Tr[M_{W'_K}^{W'_F}(\rho)(\sigma_3 a)]= Tr[Ind_{W'_{K'}}^{W'_F}(M_{W'_B}^{W'_{K'}}(\omega))(\sigma_3 a)]+Tr[ Ind_{W'_{K''}}^{W'_F}(M_{W'_B}^{W'_{K''}}(\omega))(\sigma_3 a)]=\omega(\sigma_3 a\sigma_3 a)+\omega^{\sigma_1}(\sigma_3 a\sigma_3 a)$.\\
Hence we have the isomorphism $$ M_{W'_K}^{W'_F}(\rho) \simeq
 Ind_{W'_{K'}}^{W'_F}(M_{W'_B}^{W'_{K'}}(\omega))\oplus
 Ind_{W'_{K''}}^{W'_F}(M_{W'_B}^{W'_{K''}}(\omega)).$$\\
From this we deduce that $$\boxed{L(M_{W'_K}^{W'_F}(\rho),s)=
 L(\omega_{|K'^*},s)L(\omega_{|K''^*},s).}$$

\item[ii)] Let $L$ be a quadratic extension of $F$, such that $\rho= Ind_{W'_L}^{W'_K}(\chi)$, with $\chi$ regular, is
 not isomorphic to a representation of the form
 $Ind_{W'_B}^{W'_K}(\omega)$ as in i), then
 $$\boxed{L(M_{W'_K}^{W'_F}(\rho),s)=1.}$$
 Indeed, we show that $M_{W'_K}^{W'_F}(\rho)^{I'_F}=\left\lbrace
 0\right\rbrace $. If it wasn't the case, the representation $(M_{W'_K}^{W'_F}(\rho),V)$
 would admit a $I'_F$-fixed vector, and so would its contragredient
 $V^*$. Now in the subspace of $I_F'$-fixed vectors of $V^*$, choosing an eigenvector of
 $M_{W'_K}^{W'_F}(\rho)(\phi_F)$, we would deduce the existence of a
 linear form $L$ on
 $(M_{W'_K}^{W'_F}(\rho),V)$ which transforms under $W'_F$ by an
 unramified
 character $\mu$ of $W'_F$. If we identify $\mu$ with a character $\mu'$ of $F^*$, the restriction of $\mu$ to $W'_K$ corresponds to $\mu' \circ
 N_{K/F}$ of $K^*$, so we can write it as $\theta {\theta}^{\sigma}$, where $\theta$ is a character of $W'_K$ corresponding to an extension of $\mu'$ to $K^*$. As the restriction of $M_{W'_K}^{W'_F}$ to $W'_K$ is isomorphic to
 $\rho
 \otimes \rho ^{\sigma}$, we deduce that ${\theta}^{-1} \rho \otimes
 ({\theta}^{-1}\rho) ^{\sigma}$ is $W'_K$ distinguished, that is
 ${\theta}
 \rho^{\vee} \simeq {({\theta}^{-1} \rho)} ^{\sigma}$. But from the proof of Theorem \ref{+dist}, this would imply that
 ${\theta}^{-1} \rho$ hence $\rho$, could be induced from a character
 of a biquadratic extension of $F$, which we supposed is not the
 case.

\item[iii)] Suppose $\rho=sp(2)$ acts on the space $\mathbb{C}^2$ with canonical basis $(e_1,e_2)$ by the natural action
 $\rho \left[ h, M \right](v)=M(v)$ for $h$ in $W_K$, $M$ in $SL(2,\mathbb{C})$ and $v$ in $\mathbb{C}^2$. Then the space of $M_{W'_K}^{W'_F}(\rho)$ is
 $V\otimes
 V$ and  $SL(2,\mathbb{C})$ acts on it as  $sp(2)\otimes sp(2)$. Decomposing $V\otimes V$ as the direct sum $Alt(V)\oplus Sym(V)$, we
 see that $SL(2,\mathbb{C})$ acts as $1$ on $Alt(V)$, and
 $M_{W'_K}^{W'_F}(\rho)\left[1,\left(\begin{array}{lll} x & 0 \\ 0&
 x^{-1} \\ \end{array}
 \right)\right] (e_1\otimes e_1)=x^{2} e_1\otimes e_1 $. Hence the representation of $SL(2,\mathbb{C})$ on $Sym(V)$ must be
 $sp(3)$. The Weil group $W_F$ acts as $\eta _{K/F}$ on $Alt(V)$ and trivially on
 $Sym(V)$, finally $M_{W'_K}^{W'_F}(\rho)$ is isomorphic to
 $sp(3)\oplus
 \eta _{K/F}$. Tensoring with a character $\chi$, we have $M_{W'_K}^{W'_F}(\chi
 sp(2))= \chi _{|F^*} M_{W'_K}^{W'_F}(sp(2))= \chi _{|F^*}\eta _{K/F}
 \oplus
 \chi _{|F^*} sp(3) $. Hence one has the following equality: $$\boxed{L(M_{W'_K}^{W'_F}(\chi
 sp(2)),s)=L(\chi_{|F^*} \eta _{K/F},s)L(\chi_{|F^*},s+1).}$$

\item[iv)] If $\rho=\lambda \oplus \mu$, with $\lambda$ and $\mu$ two
 characters of $W'_K$, then from \cite{P}, Lemma 7.1, we have
 $M_{W'_K}^{W'_F}(\rho)=\lambda _{|F^*} \oplus {\mu} _{|F^*} \oplus
 Ind_{W'_K}^{W'_F}(\lambda{\mu}^{\sigma})$. Hence we have

 $$\boxed{L(M_{W'_K}^{W'_F}(\rho))=L(\lambda_{|F^*},s)L(\mu_{|F^*},s)L(\lambda{\mu}^{\sigma},s).}$$  

 \end{description}

\subsection{Asai $L$-functions for ordinary representations of
 $GL(2)$} 

In this subsection, we compute Asai $L$-functions for ordinary (i.e. non exceptional) representations of $G_2(K)$, and prove (Theorem \ref{egal}) that they are equal to the corresponding functions $L_W$ of imprimitive representations of $W'_K$.\\

In order to compute $L_{As}$, we first compute $L_1$, but this latter computation is easy because Kirillov models of infinite-dimensional irreducible representations of $G_2(K)$ are well-known (see \cite{Bu}, Th. 4.7.2 and 4.7.3).\\
Let $\pi$ be an irreducible infinite-dimensional (hence generic) representation of $G_2(K)$, we have the following situations for the computation of $L_1(\pi,s)$.

\begin{description}
 
 \item[i) and ii)] If $\pi$ is supercuspidal, its Kirillov model consists of functions with compact support on $K^*$, hence $$L_1(\pi,s)=1.$$

\item[iii)] If $\pi=\sigma(\chi)$ ($\sigma(\chi{|\ |_K}^{1/2}, \chi {|\ |_K}^{-1/2})$ in \cite{Bu}) is a special series representation of $G_2(K)$, twist of the Steinberg representation by the character $\chi$ of $K^*$, the Kirillov model of $\pi$ consists of functions of $D(K)$ multiplied by $\chi{|\ |_K}$. Hence their restrictions to $F$ are functions of $D(F)$ multiplied by $\chi{|\ |_F}^2$, and the ideal $I_1 (\pi)$ is generated by functions of $s$ of the form $$\int_{F^*}\phi(t)\chi(t){|t|_F}^{s-1}{|t|_F}^2 d^*t= \int_{F^*}\phi(t)\chi(t){|t|_F}^{s+1}d^*t,$$ for $\phi$ in $D(F)$, hence we have $$L_1(\pi,s)=L(\chi_{|F^*},s+1).$$ 

\item[iv)] If $\pi=\pi(\lambda, \mu)$ is the principal series representation ($\lambda$ and $\mu$ being two characters of $K^*$, with $\lambda \mu^{-1}$ different from $|\ |$ and $|\ |^{-1}$) corresponding to the representation $\lambda \oplus \mu$ of $W'_K$.\\

If $\lambda \neq \mu$, the Kirillov model of $\pi$ is given by functions of the form ${|\ |_K}^{1/2}\chi \phi_1 + {|\ |_K}^{1/2}\mu \phi_2$, for $\phi_1$ and $\phi_2$ in $D(K)$, and $$L_1(\pi,s)=L(\lambda_{|F^*},s)\vee L(\mu_{|F^*},s).$$

If $\lambda = \mu$, the Kirillov model of $\pi$ is given by functions of the form ${|\ |_K}^{1/2}\lambda \phi_1 + {|\ |_K}^{1/2}\lambda v_K (t) \phi_2$, for $\phi_1$ and $\phi_2$ in $D(K)$, and $$L_1(\pi,s)=L(\lambda_{|F^*},s)^2.$$
\end{description}

In order to compute $L_{rad(ex)}$ for ordinary representations, we need to know when they are distinguished by a character $|\ |_{F}^{- s_0}$ for some $s_0$ in $\mathbb{C}$, we will then use Theorem \ref{s-dist}. The answer is given by the following, which is a mix of Theorem \ref{distcusp} and Proposition B.17 of \cite{FH}:

\begin{thm}\label{cardist}
\begin{description}
 \item[a)] A dihedral supercuspidal representation $\pi$ of $G_2(K)$ is $|\ |_{F}^{- s_0}$-distinguished if and only if there exists a quadratic extension $B$ of $K$, biquadratic over $F$ (hence there are two other extensions between $F$ and $B$ that we call $K'$ and $K''$), and a character of $B^*$ regular with respect to $N_{B/K}$ which restricts either to $K'$ as $|\ |_{K'}^{-s_0}$ or to $K''$ as $|\ |_{K''}^{-s_0}$, such that $\pi$ is equal to $\pi(\omega)$.
\item[b)] Let $\mu$ be a character of $K^*$, then the special series representation $\sigma(\mu)$ is $|\ |_{F}^{- s_0}$-distinguished if and only if $\mu$ restricts to $F^*$ as $\eta_{K/F}|\ |_{F}^{- s_0}$.
\item[c)] Let $\lambda$ and $\mu$ be two characters of $K^*$, with $\lambda\mu^{-1}$ and $\lambda^{-1}\mu$ different from $|\ |_{K}$, then the principal series representation $\pi(\lambda,\mu)$ is $|\ |_{F}^{- s_0}$-distinguished if and only if either $\lambda$ and $\mu$ restrict as $|\ |_{F}^{- s_0}$ to $F^*$ or $\lambda\mu^{\sigma}$ is equal to $|\ |_{K}^{- s_0}$.
\end{description}
\end{thm}
\begin{proof} Let $\pi$ be a representation, it is $|\ |_{F}^{- s_0}$-distinguished if and only if $\pi \otimes |\ |_{K}^{ s_0/2}$ is distinguished because $|\ |_{K}^{- s_0/2}$ extends $|\ |_{F}^{- s_0}$, it then suffices to apply Theorem \ref{distcusp} and Proposition B.17 of \cite{FH}. We give the full proof for case a). Suppose $\pi$ is dihedral supercuspidal and $\pi \otimes |\ |_{K}^{ s_0/2}$ is distinguished. From Theorem \ref{distcusp}, the representation $\pi \otimes |\ |_{K}^{ s_0/2}$ must be of the form $\pi(\omega)$, for $\omega$ a character of quadratic extension $B$ of $K$, biquadratic over $F$, such that if we call $K'$ and $K''$ two other extensions between $F$ and $B$, $\omega$ doesn't factorize through $N_{B/K}$ and restricts either trivially on $K'^*$, or trivially on $K''^*$. But $\pi$ is equal to $\pi(\omega)\otimes |\ |_{K}^{-s_0/2}= \pi(\omega |\ |_{B}^{-s_0/2})$ because $|\ |_{B}=|\ |_{K} \circ N_{B/K}$. As $|\ |_{B}^{-s_0/2}$ restricts to $K'$ (resp. $K''$) as $|\ |_{K'}^{-s_0}$ (resp. $|\ |_{K''}^{-s_0}$), case a) follows.\end{proof}

We are now able to compute $L_{rad(ex)}$, hence $L_{As}$ for ordinary representations.

\begin{description}

\item[i)] Suppose that $\pi= \pi(Ind_{W'_B}^{W'_K}(\omega))=\pi(\omega)$ is supercuspidal, with Langlands parameter $Ind_{W'_B}^{W'_K}(\omega)$, where $\omega$ is a multiplicative character of a biquadratic extension $B$ over $F$ that doesn't factorize through $N_{B/K}$.\\
We denote by $K'$ and $K''$ the two other extensions between $B$ and $F$. Here $L_1(\pi,s)$ is equal to one.\\
We have the following series of equivalences:
 $$\aligned s_0 \ is \ a \ pole \ of \ L_{As}(\pi(\omega),s)&\Longleftrightarrow \pi(\omega) \ is \ |\ |_F^{-s_0} -distinguished \\
                                                                &\Longleftrightarrow \omega_{|K'^*} =|\ |_{K'}^{-s_0} \ or \ \omega_{|K''^*}=|\ |_{K''}^{-s_0} \\
                             &\Longleftrightarrow  s_0 \ is \ a \ pole \ of \ L(\omega_{|K'^*},s) \ or \ of \ L(\omega_{|K''^*},s)\\
                                    &\Longleftrightarrow  s_0 \ is \ a \ pole \ of \ L(\omega_{|K'^*},s)\vee L(\omega_{|K''^*},s)
\endaligned$$

As both functions $L_{As}(\pi(\omega),s)$ and $L(\omega_{|K'^*},s)\vee L(\omega_{|K''^*},s)$ have simple poles and are Euler factors, they are equal. Now suppose that $L(\omega_{|K'^*},s)$ and $L(\omega_{|K''^*},s)$ have a common pole $s_0$, this would imply that $\omega_{|K'^*} =|\ |_{K'}^{-s_0}$ and $\omega_{|K''^*}=|\ |_{K''}^{-s_0}$, which would mean that $\omega|\ |_B^{s_0 /2}$ is trivial on $K'^*K''^*$. According to Lemma \ref{Ker}, this would contradict the fact that $\omega$ does not factorize through $N_{B/K}$, hence $L(\omega_{|K'^*},s)\vee L(\omega_{|K''^*},s)=L(\omega_{|K'^*},s)L(\omega_{|K''^*},s)$. Finally we proved:
$$\boxed{L_{As}(\pi(\omega),s)=L(\omega_{|K'^*},s)L(\omega_{|K''^*},s).}$$

\item[ii)] Suppose that $\pi$ is a supercuspidal representation, corresponding to an imprimitive representation of $W'_K$ that cannot be induced from a character of the Weil-Deligne group of a biquadratic extension of $F$.
Then necessarily $\pi$ cannot be $|\ |_F^{-s_0}$-distinguished, for any complex number $s_0$ of $\mathbb{C}$.\\
If it was the case, from Theorem \ref{cardist}, it would correspond to a Weil representation $\pi(\omega)$ for some multiplicative character of a biquadratic extension of $F$, which cannot be. Hence $L_{rad(ex)}(\pi,s)$ has no pole and is equal to one because it is an Euler factor, so we proved that:  
 $$\boxed{L_{As}(\pi,s)=1.}$$ 

\item[iii)] If $\pi$ is equal to $\sigma(\chi)$, then $L_1(\pi,s)=L(\chi_{|F^*},s+1)$. We want to compute $L_{rad(ex)}(\pi,s)$, we have the following series of equivalences:
$$\aligned s_0 \ is \ an \ exceptional \ pole \ of \ L_{As}(\sigma( \chi),s) &\Longleftrightarrow \sigma( \chi) \ is \ |\ |_F^{-s_0} -distinguished \\
&\Longleftrightarrow \chi_{|F^*}=\eta_{K/F}|\ |_F^{-s_0}\\
&\Longleftrightarrow s_0 \ is \ a \ pole \ of \ L(\chi_{|F^*}\eta_{K/F},s)
\endaligned$$

As both functions $L_{rad(ex)}(\pi,s)$ and $L(\chi_{|F^*}\eta_{K/F},s)$ have simple poles and are Euler factors, they are equal, we thus have:
$$\boxed{L_{As}(\sigma( \chi)=L(\chi_{|F^*},s+1)L(\chi_{|F^*}\eta_{K/F},s).}$$

\item[iv)] If $\pi=\pi(\lambda, \mu)$, we first compute $L_{rad(ex)}(\pi,s)$. We have the following series of equivalences:
$$\aligned s_0 \ is \ an \ exceptional \ pole \ of \ L_{As}(\pi(\lambda, \mu),s)&\Longleftrightarrow  \pi(\lambda, \mu)  \ is \ |\ |_F^{-s_0} -distinguished \\
 &\Longleftrightarrow \lambda\mu^{\sigma}=|\ |_K^{-s_0} \ or, \ \lambda_{|F^*}=|\ |_F^{-s_0} \ and \ \mu_{|F^*}=|\ |_F^{-s_0}\\
  &\Longleftrightarrow s_0 \ is \ a \ pole \ of \ L(\lambda\mu^{\sigma},s) \ or \ of \ L(\lambda_{|F^*},s)\wedge L(\mu_{|F^*},s)\\
      &\Longleftrightarrow s_0 \ is \ a \ pole \ of \ L(\lambda\mu^{\sigma},s)\vee[L(\lambda_{|F^*},s)\wedge L(\mu_{|F^*},s)]
\endaligned$$

As both functions $L_{rad(ex)}(\pi(\lambda, \mu),s)$ and $L(\lambda\mu^{\sigma},s)\vee[L(\lambda_{|F^*},s)\wedge L(\mu_{|F^*},s)]$ have simple poles and are Euler factors, they are equal.

If $\lambda \neq \mu$, then $L_1(\pi,s)=L(\lambda_{|F^*},s)\vee L(\mu_{|F^*},s)$. But $L(\lambda\mu^{\sigma},s)$ and $L(\lambda_{|F^*},s)\wedge L(\mu_{|F^*},s)$ have no common pole. If there was a common pole $s_0$, one would have $\lambda\mu^{\sigma}=|\ |_K^{-s_0}$, $\lambda_{|F^*}=|\ |_F^{-s_0}$ and $\mu_{|F^*}=|\ |_F^{-s_0}$. From $\mu_{|F^*}=|\ |_F^{-s_0}$, we would deduce that $\mu \circ N_{K/F}= |\ |_K^{-s_0}$, i.e. $\mu^{\sigma}= |\ |_K^{-s_0} \mu^{-1}$, and $\lambda\mu^{\sigma}=|\ |_K^{-s_0}$ would imply $\lambda = \mu$, which is absurd. Hence $L_{rad(ex)}(\pi,s)=L(\lambda\mu^{\sigma},s)[L(\lambda_{|F^*},s)\wedge L(\mu_{|F^*},s)]$, and finally we have $L_{As}(\pi,s)=L_1(\pi,s)L_{rad(ex)}(\pi,s)=L(\lambda_{|F^*},s)L(\mu_{|F^*},s)L(\lambda{\mu}^{\sigma},s)$.\\
If $\lambda$ is equal to $\mu$, then $L_1(\pi,s)=L(\lambda_{|F^*},s)^2$, and $L_{rad(ex)}(\pi(\lambda, \mu),s)=L(\lambda \circ N_{K/F},s)\vee L(\lambda_{|F^*},s)$. As $L(\lambda \circ N_{K/F},s)=L(\lambda_{|F^*},s)L(\eta_{K/F}\lambda_{|F^*},s)$, we have $L_{rad(ex)}(\pi(\lambda, \mu),s)=L(\lambda \circ N_{K/F},s)$. Again we have $L_{As}(\pi,s)=L(\lambda_{|F^*},s)L(\mu_{|F^*},s)L(\lambda{\mu}^{\sigma},s)$.\\
In both cases, we have $$\boxed{L_{As}(\pi(\lambda, \mu),s)= L(\lambda_{|F^*},s)L(\mu_{|F^*},s)L(\lambda{\mu}^{\sigma},s).}$$

\end{description}

Eventually, comparing with equalities of subsection \ref{Weil}, we proved the following:

\begin{thm}\label{egal}
 
Let $\rho \mapsto \pi(\rho)$ be the Langlands correspondence from two dimensional representations of $W'_K$ to smooth irreducible infinite-dimensional representations of $G_2(K)$, then if $\rho$ is not primitive, $\pi(\rho)$ is ordinary and we have the following equality of $L$-functions:

$$\boxed{L_{As}(\pi(\rho),s)=L(M_{W'_K}^{W'_F}(\rho),s)}$$

\end{thm}

As said in the introduction, combining Theorem 1.6 of \cite{AR} and Theorem of pargraph 1.5 in \cite{He}, one gets that $L(M_{W'_K}^{W'_F}(\rho),s)=L_{As}(\pi(\rho),s)$ for $\pi(\rho)$ a discrete series representation, so that we have actually the following:

\begin{thm}
 
Let $\rho \mapsto \pi(\rho)$ be the Langlands correspondence from two dimensional representations of $W'_K$ to smooth irreducible infinite-dimensional representations of $G_2(K)$, we have the following equality of $L$-functions:

$$\boxed{L_{As}(\pi(\rho),s)=L(M_{W'_K}^{W'_F}(\rho),s)}$$\end{thm}

\vskip3em

\begin{concl}

The results of Section \ref{gl(2)} give a local proof of the equality of $L_W$ and $L_{As}$, and effective computations of these functions. As it was said in the introduction, the latter equality is known for discrete series representations of $G_n(K)$ but the proof is global. Hence the essentially new information is the equality for principal series representations of $G_2(K)$.\\
Now the following conjecture is expected to be true:

\begin{conj}\label{distgen} Let $(n_1,\dots,n_t)$ be a partition of $n$, and for each $i$ between $1$ and $t$, let $\Delta_i$ be a quasi-square-integrable representation of $G_{n_i}(K)$. The generic representation $\pi$ of $G_n(K)$ obtained by normalised parabolic induction of the $\Delta_i$'s is distinguished if and only if there is a reordering of these representations and an integer $r$ between $1$ and $t/2$, such that $\Delta_{i+1}^{\sigma} = \Delta_i^{\vee} $ for $i=1,3,..,2r-1$, and $\Delta_{i}$ is distiunguished for $i > 2r$.\end{conj}

In a work to follow, we intend to prove that assuming this conjecture, the functions $L_W$ and $L_{As}$ agree on generic representations of $G_n(K)$. As Conjecture \ref{distgen} is proved in \cite{M} for principal series representations, this would give the equality of the $L$ functions for  principal series representations of $G_n(K)$.  
\end{concl}

\section{Appendix. Dihedral supercuspidal distinguished representations}\label{appendix}

The aim of this section is to give a description of dihedral supercuspidal distinguished representations of $G_2(K)$ in terms of Langlads parameter, it is done in Theorem \ref{distcusp}. 

\subsection{Preliminary results}

Let $E$ be a local field, $E'$ be a quadratic extension of $E$, $\chi$ a character of $E^*$, $\pi$ be a smooth irreducible infinite-dimensional representation of $G_2(E)$, and $\psi$ a non trivial character of $E$.\\
 We denote by $L(\chi,s)$ and $\epsilon(\chi,s,\psi)$ the functions of the complex variable $s$ defined before Proposition 3.5 in \cite{JL}.We denote by $\gamma(\chi,s,\psi)$ the ratio $\epsilon(\chi,s,\psi)L(\chi,s)/L(\chi^{-1},1-s)$.\\
 We denote by $L(\pi,s)$ and $\epsilon(\pi,s,\psi)$ the functions of the complex variable $s$ defined in Theorem 2.18 of \cite{JL}. We denote by $\gamma(\pi,s,\psi)$ the ratio $\epsilon(\pi,s,\psi)L(\pi,s)/L(\pi^{\vee},1-s)$.\\
We denote by $\lambda(E'/E,\psi)$ the Langlands-Deligne factor defined before Proposition 1.3 in \cite{JL}, it is equal to $\epsilon(\eta_{E'/E},1/2,\psi)$. As $\eta_{E'/E}$ is equal to $\eta_{E'/E}^{-1}$, the factor $\lambda(E'/E,\psi)$ is also equal to $\gamma(\eta_{E'/E},1/2,\psi)$.\\
From Theorem 4.7 of \cite{JL}, if $\omega$ is a character of $E'^*$, then $L(\pi(\omega),s)$ is equal to $L(\omega,s)$, and $\epsilon(\pi,s,\psi)$ is equal to $\lambda(E'/E,\psi)\epsilon(\pi,s,\psi)$, hence $\gamma(\pi,s,\psi)$ is equal to $\lambda(E'/E,\psi)\gamma(\pi,s,\psi)$.\\

We will need four results. The first is due to Fröhlich and Queyrut, see \cite{D} Theorem 3.2 for a
 quick proof using a Poisson formula:

\begin{prop}\label{Fro}
Let $E$ be a local field, $E'$ be a quadratic extension of $E$, $\chi'$ a character of $E'^*$ trivial on $E^*$, and $\psi'$ a non trivial character of $E'$ trivial on $E$, then $\gamma(\chi',1/2 ,\psi')=1$.
\end{prop}

The second is a criterion of Hakim:

\begin{thm}\label{Hakim}{(\cite{Ha}, Theorem 4.1)}
Let $\pi$ be an irreducible supercuspidal representation of $G_2(K)$ with
 central character trivial on $F^*$, and $\psi$ a nontrivial character
 of $K$ trivial on $F$. Then  $\pi$ is distinguished if and only if
 $\gamma (\pi \otimes \chi ,1/2 ,\psi)=1 $ for every character $\chi$ of $K^*$ trivial
 on $F^*$. 
\end{thm}

The third is due to Flicker:
\begin{thm}\label{Flic}{(\cite{F1}, Proposition 12)}
\label{autodual}
 Let $\pi$ be a smooth irreducible distinguished representation of $G_n(K)$, then $\pi
 ^{\sigma}$ is isomorphic to $\pi^{\vee}$.
\end{thm}

The fourth is due to Kable in the case of $G_n(K)$, see \cite{AT} for a local proof in the case of $G_2(K)$:
\begin{thm}\label{exclusive}{(\cite{AT}, Proposition 3.1}
There exists no supercuspidal representation of $G_2(K)$ which is
 distinguished and $\eta _{K/F}$-distinguished at the same time.
\end{thm}

 \subsection{Distinction criterion for dihedral supercuspidal representations}

As a dihedral representation's parameter is a multiplicative character of a quadratic extension $L$ of $K$, we first look at the properties of the tower $F\subset K\subset L$.
Three cases arise:

\begin{enumerate}

 \item  $L /F$ is biquadratic (hence Galois), it contains $K$ and two
 other quadratic extensions $F$, $K'$ and $K''$. 

\begin{figure}[here]\label{fig1}
 \centering
\includegraphics[width=2.1in,height=1in]{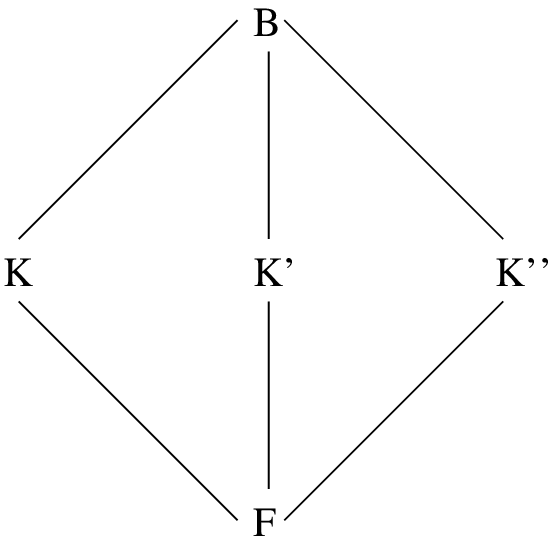}
\caption{}
\end{figure}

Its Galois group is isomorphic with $\mathbb{Z} / 2 \mathbb{Z} \times
 \mathbb{Z} / 2 \mathbb{Z}$, its non trivial elements are $\sigma_{L/K}$, $\sigma_{L/K'}$ and $\sigma_{L/K''}$.
The conjugation $\sigma_{L/K}$ extend $\sigma_{K'/F}$ and $\sigma_{K''/F}$.

\item $L/F$ is cyclic with Galois group isomorphic with $\mathbb{Z} /
 4 \mathbb{Z}$, in this case we fix an element $\tilde{\sigma}$ in $G(L/F)$ extending $\sigma$, it is of order 4.

 \item $L/F$ non Galois. Then its Galois closure $M$ is quadratic over
 $L$ and the Galois group of $M$ over $F$ is dihedral with order 8.
 To see this, we consider a morphism $\tilde{\theta}$ from $L$ to
 $\bar{F}$ which extends $\theta$. Then if $L' = \tilde{\theta}(L)$,
 $L$ and $L'$ are distinct, quadratic over $K$ and generate $M$ biquadratic
 over $K$. $M$ is the Galois closure of $L$ because any morphism from $L$
 into $\bar{F}$, either extends $\theta$, or the identity map of $K$,
 so that its image is either $L$ or $L'$, so it is always included in
 $M$. Finally the Galois group $M$ over $F$ cannot be abelian (for $L$ is not
 Galois over $F$), it is of order $8$, and it's not the quaternion group which
 only has one element of order $2$, whereas here $\sigma_{M/L}$ and $\sigma_{M/L'}$ are of order $2$. Hence it is the dihedral group of
 order 8 and we have the folowing lattice, where $M/K'$ is cyclic of degree 4, $M/K$ and $B/F$ are biquadratic.\\

\begin{figure}[here]\label{fig2}
 \centering
\includegraphics[width=3in,height=1.2in]{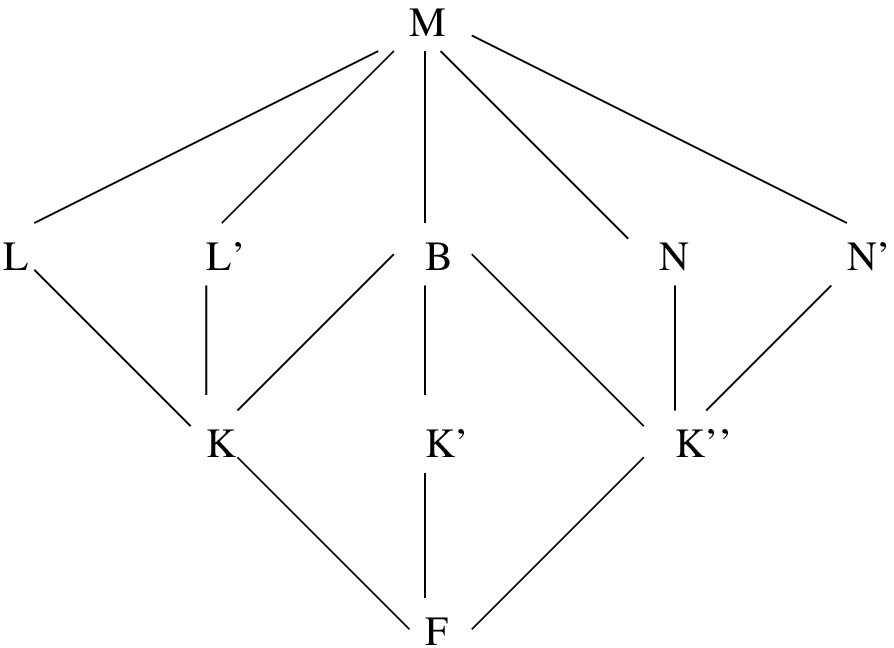}
\caption{}
\end{figure} 
\end{enumerate}

We now prove the following proposition:

\begin{prop}\label{+dist}

 If a supercuspidal dihedral representation $\pi$ of $G_2(K)$
 verifies $\pi^{\vee}={\pi}^{\sigma}$, there exists a biquadratic
 extension $B$ of $F$, containing $K$, such that if we call $K'$ and $K''$ the two other extensions between $F$ and $B$, there is a character $\omega$ of $B$ trivial either on $N
 _{B/K'}({B}^*)$ or on $N_{B/K''}({B}^*)$, such that $\pi = \pi
 (\omega)$.

\end{prop}

\begin{proof}Let $L$ be a quadratic extension of $K$ and $\omega$ a regular
 multiplicative of $L$ such that $\pi = \pi (\omega)$, we denote by $\sigma$
 the
 conjugation of $L$ over $K$, three cases show up:

\begin{enumerate}

\item $L/F$ is biquadratic.
The conjugations $\sigma_{L/K'}$ and $\sigma_{L/K''}$ both extend
 $\sigma$, hence from Theorem 1 of \cite{GL}, we have $\pi (\omega)^{\sigma}= \pi (\omega^{\sigma_{L/K'}})$.
 The condition ${\pi}^{\vee} = {\pi}^{\sigma}$ which one can also read
 $\pi ({\omega}^{-1})= \pi ({\omega}^{\sigma_{L/K'}})$, is then
 equivalent from Appendix B, (2)b)1) of \cite{GL}, to ${\omega}^{\sigma_{L/K'}}= {\omega}^{-1}$ or
 ${\omega}^{\sigma_{L/K''}}={\omega}^{-1}$.
 This is equivalent to $\omega$ trivial on $N _{L/K'}({L}^*)$ or on $ N
 _{L/K''}({L}^*)$.

\item $L/F$ is cyclic, the regularity of $\omega$ makes the condition
 $\pi ({\omega}^{-1}) = \pi ({\omega})^{\sigma}$
 impossible.
Indeed one would have from Theorem 1 of \cite{GL} $\pi({\omega}^{\tilde{\sigma}})=\pi({\omega}^{-1})$, which from Appendix B, (2)b)1) of \cite{GL} would imply ${\omega}^{\tilde{\sigma}}= \omega$ or ${\omega}^{\tilde{\sigma}^{-1}}= \omega$. As ${\tilde{\sigma}}^2 = {\tilde{\sigma}}^{-2}=\sigma$, this would in turn imply $\omega^{\sigma}=\omega$, and $\omega$ would be trivial on the kernel of $N _{L/K}$ according to Hilbert's Theorem 90.
$\pi^{\vee}$ can therefore not be isomorphic to ${\pi}^{\sigma}$.

\item $L/K$ is not Galois (which implies $q \equiv 3[4]$ in the case p
 odd).
 Let ${\pi}_{B/K}$ be the representation of ${G_2 (B)}$ which is the
 base change lift of ${\pi}$ to $B$. As ${\pi}_{B/K}= \pi ({\omega}
 \circ N _{M/L})$, if ${\omega} \circ N _{M/L} = {\mu}\circ N _{M/B} $
 for a character $\mu$ of $B^*$, then $\pi(\omega) = \pi(\mu)$
 (cf.\cite{GL}, (3) of Appendix B) and we are brought back to case 1.
 Otherwise ${\omega} \circ N _{M/L}$ is regular with respect to $N _{M/B} $.
 If  ${\pi}^{\sigma} ={\pi}^{\vee}$, we would
 have ${\pi}_{B/K}^{\sigma_{B/K'}}= {\pi}_{B/K}^{\vee}$ from Theorem 1 of
 \cite{GL}.
 That would contradict case 2 because $M/K'$ is cyclic.\end{enumerate} \end{proof}

 We described in the previous proposition representations $\pi$ of
 $G_2(K)$ verifying $\pi^{\vee}={\pi}^{\sigma}$, now we characterize those
 who are ${G _2 (F)}$-distinguished among
 them (from Theorem \ref{Flic}, a distinguished representation
 always satisfies the previous condition).\\

\begin{thm}\label{distcusp}
A dihedral supercuspidal representation $\pi$ of $G_2 (K)$ is $G_2
 (F)$-distinguished if and only if there exists a quadratic extension
 $B$ of $K$ biquadratic over $F$ such that if we call $K'$ and $K''$ the two other extensions between $B$ and $F$, there is character $\omega$ of $B^*$ that does not factorize through $N_{B/K}$ and trivial either on ${K'}^*$ or on ${K''}^*$, such that $\pi = \pi (\omega)$.
\end{thm}

\begin{proof} From Theorem \ref{Flic} and Proposition \ref{+dist}, we can suppose that $\pi
 = \pi (\omega)$, for
 $\omega$ a regular multiplicative character of a quadratic extension
 $B$ of $K$ biquadratic over $F$, with $\omega$ trivial on $N
 _{L/K'}({K'}^*)$
 or on $ N _{B/K''}({K''}^*)$. We will need the following:
\begin{LM}\label{normbiquad}
 Let $B$ be a quadratic extension of $K$ biquadratic over $F$, then $F^*$ is a subset of $N_{B/K}(B^*)$
\end{LM}
\begin{proof}[Proof of Lemma \ref{normbiquad}]
The group $N_{B/K}(B^*)$ contains the two groups $N_{B/K}(K'^*)$ and $N_{B/K}(K''^*)$, which, as $\sigma_{B/K}$ extends $\sigma_{K'/F}$ and $\sigma_{K''/F}$, are respectively equal to $N_{K'/F}(K'^*)$ and $N_{K''/F}(K''^*)$. But these two groups are distinct of index 2 in $F^*$ from local cassfield theory, thus they generate $F^*$, which is therefore contained in $N_{B/K}(B^*)$. 
\end{proof}

 We now choose $\psi$ a non trivial character of $K/F$ and denote by $\psi _B$
 the character $\psi \circ Tr_{B/K}$, it is trivial on $K'$ and $K''$.\\
 Suppose $\omega$ trivial on $K'$ or $K''$, then the restriction of
 the central character $\eta_{B/K} \omega$ of $\pi(\omega)$ is trivial
 on $F^*$ according to Lemma \ref{normbiquad}.\\
As we have $\gamma(\pi(\omega),1/2 ,\psi)=\lambda(B/K,\psi)\gamma(\omega,1/2 ,\psi_B)=\gamma(\eta_{B/K},1/2 ,\psi)\gamma(\omega,1/2 , \psi_B)$, we deduce from Lemma \ref{normbiquad} and Proposition \ref{Fro} that $\gamma(\pi(\omega),1/2 ,\psi)$ is equal to one, hence from Theorem \ref{Hakim}, the representation $\pi(\omega)$ is distinguished.\\
Now suppose $\omega|_{K'}=\eta_{B/K'}$ or $\omega|_{K''}=\eta_{B/K''}$, let $\chi$ be a character of $K^*$ extending $\eta_{K/F}$, then $\pi(\omega) \otimes \chi = \pi(\omega \chi \circ N_{B/K})$. As ${N_{B/K}}_{|K'}=N_{K'/F}$ and ${N_{B/K}}_{|K''}=N_{K''/F}$, we have $\chi \circ {N_{B/K}}_{|K'}=\eta_{B/K'}$ and $\chi \circ {N_{B/K}}_{|K''}=\eta_{B/K''}$, hence from what we've just seen, $\pi(\omega) \otimes \chi$ is distinguished, i.e. $\pi(\omega)$ is $\eta_{K/F}$-distinguished.\\
From Theorem \ref{exclusive}, $\pi$ cannot be distinguished and $\eta_{K/F}$-distinguished at the same time, and the theorem follows. \end{proof}

We end with the following lemma:
\begin{LM}\label{Ker}
Let $B$ be a quadratic extension of $K$ which is biquadratic over $F$. Call $K'$ and $K''$ the two other extensions between $F$ and $B$, then the kernel of $N_{B/K}$ is a subgroup of the group $N_{B/K'}(B^*)N_{B/K''}(B^*)$. 
\end{LM}
\begin{proof} If $u$ belongs to $Ker(N_{B/K})$, it can be written $x/\sigma_{B/K}(x)$ for some $x$ in $B^*$ according to Hilbert's Theorem 90. Hence we have $u=(x\sigma_{B/K'}(x))/(\sigma_{B/K}(x)\sigma_{B/K'}(x))=N_{B/K'}(x)/N_{B/K''}(\sigma_{B/K}(x))$, and $u$ belongs to $N_{B/K'}(B^*)N_{B/K''}(B^*)$. \end{proof}

\begin{cor}
The (either/or) in Proposition \ref{+dist} and Theorem \ref{distcusp} is exclusive
\end{cor}
\begin{proof} In fact, in the situation of Lemma \ref{Ker}, a character $\omega$ that is trivial on $N_{B/K'}(B^*)$ and $N_{B/K''}(B^*)$ factorizes through $N_{B/K}$, and $\pi(\omega)$ is not supercuspidal.\end{proof}

\section*{Acknowledgements}
I would like to thank Corinne Blondel and Paul G\'{e}rardin for many helpful comments. I also thank Jeffrey Hakim who allowed me to have access to Youngbin Ok's PHD thesis.

\end{document}